\documentclass{amsart}
\newtheorem{lem}{Lemma}[section]
\newtheorem{thm}{Theorem}[section]
\newtheorem{prop}{Proposition}[section]

\newtheorem{cor}{Corollary}[section]

\def\cit{{\mathbb C}}

\title[B-Fredholm and Drazin invertible operators through localized SVEP]
{B-Fredholm and Drazin invertible operators through localized SVEP}
\author{M. Amouch and H. Zguitti}
\address{M. Amouch: Department of Mathematics, Faculty of science
Semlalia, B.O : 2390 Marrakesh, Morocco.}
\email{m.amouch@ucam.ac.ma}
\address{H. Zguitti: Department of Mathematics, Faculty of science of Rabat,
B.O : 1014 Rabat, Morocco.} \email{zguitti@hotmail.com}

\begin{document}


\maketitle \thispagestyle{empty}

\begin{abstract}
Let $X$ a Banach space and $T$ a bounded linear operator on $X.$
We denote by $S(T)$ the set of all $\lambda \in \cit$ such that
$T$ does not have the single-valued extension property at
$\lambda$. In this note we prove equality up to $S(T)$ between the
left Drazin spectrum and the left B-Fredholm spectrum and between
the semi-essential approximate point spectrum and the left Drazin
spectrum. As applications we investigate generalized Weyl's
theorem for operator matrices and multipliers operators.
\end{abstract}

\section{Introduction}

Throughout this paper, $X$ and $Y$ are Banach spaces and let ${\mathcal B}(X,Y)$ denote
 the space of all bounded linear operators from $X$ to $Y$. For $Y=X$ we write
 ${\mathcal B}(X,Y)={\mathcal B}(X)$. For $T\in {\mathcal B}(X)$, let $T^*,\, N(T)$, $R(T)$,
  $\sigma(T)$, $\sigma_s(T)$,  $\sigma_p(T)$ and $\sigma_a(T)$ denote the adjoint, the null space, the range, the spectrum, the surjective spectrum, the point spectrum and the approximate point spectrum of $T$ respectively. Let $\alpha(T)$ and $\beta(T)$ be the nullity
and the deficiency of $T$ defined by $\alpha(T)=\mbox{dim}
N(T),\mbox{ and }\beta(T)=\mbox{codim} R(T).$ If the range $R(T)$
of $T$ is closed and $\alpha(T)<\infty$ (resp. $\beta(T)<\infty$),
then $T$ is called an  {\it upper} (resp. a  {\it lower})  {\it semi-Fredholm}
operator. In the sequel $SF_+(X)$ (resp. $SF_-(X)$) will denote
the set of all upper (resp. lower ) semi-Fredholm operators. If
$T\in {\mathcal B}(X)$ is either upper or lower semi-Fredholm,
then $T$ is called a  {\it semi-Fredholm} operator, and the  {\it index} of
$T$ is defined by $\mbox{ind}(T)=\alpha(T)-\beta(T).$ If both
$\alpha(T)$ and $\beta(T)$ are finite, then $T$ is a called a
{\it Fredholm} operator. An operator $T$ is called  {\it Weyl} if it is
Fredholm
 of index zero. The Weyl spectrum $\sigma_w(T) $is defined by $
\sigma_w(T)=\{\lambda\in\mathbb{C}:T-\lambda I \mbox{ is not
Weyl}\}$.

For $T\in {\mathcal{B}}(X)$ and a nonnegative integer $n$ define
$T_{[n]}$ to be the restriction of $T$ to $R(T^n)$ viewed as a map
from $R(T^n)$ into $R(T^n)$ ( in particular $T_{[0]}=T$). If for
some integer $n$ the range space $R(T^n)$ is closed and $T_{[n]}$
is an upper (resp. a lower) semi-Fredholm operator, then $T$ is
called an {\it upper} (resp. a {\it lower}) {\it semi-B-Fredholm}
operator. In this case the {\it index} of $T$ is defined to be the
index of the semi-Fredholm operator $T_{[n]}.$ Moreover if
$T_{[n]}$ is a Fredholm operator, then $T$ is called a {\it
B-Fredholm} operator. A {\it semi-B-Fredholm} operator is an upper
or a lower semi-B-Fredholm operator (\cite{Be4, Be1, BK}. The {\it upper
semi-B-Fredholm spectrum} $\sigma_{lBF}(T)$, {\it lower
semi-B-Fredholm spectrum} $\sigma_{rBF}(T)$ and the {\it
B-Fredholm spectrum} $\sigma_{BF}(T)$ of $T$ are defined by
$$\sigma_{lBF}(T)=\{\lambda\in\cit\,:\,T-\lambda I\mbox{ is not a
semi-B-Fredholm operator}\},$$
$$\sigma_{rBF}(T)=\{\lambda\in\cit\,:\,T-\lambda I\mbox{ is not a
semi-B-Fredholm operator}\},$$
$$\sigma_{BF}(T)=\{\lambda\in\cit\,:\,T-\lambda I\mbox{ is not a
B-Fredholm operator }\}.$$ We have
$$\sigma_{BF}(T)=\sigma_{lBF}(T)\cup\sigma_{rBF}(T).$$

An operator $T\in{\mathcal L}(X)$ is said to be a {\it B-Weyl}
operator if it is a B-Fredholm operator of index zero. The {\it
B-Weyl spectrum} $\sigma_{BW}(T)$ of $T$ is defined by
$$\sigma_{BW}(T)=\{\lambda\in\cit\,:\,T-\lambda\mbox{ is not a
B-Weyl operator}\}.$$
From \cite[Lemma 4.1]{Be1}, $T$ is a B-Weyl operator if and only if $T=F\oplus N$, where $F$ is a
Fredholm operator of index zero and $N$ is a nilpotent operator.

We shall denote by $SBF_{+}^{-}(X)$
(resp. $SBF_{-}^{+}(X)$)  the class of all $T$ upper semi-B-Fredholm operators
(resp. $T$ lower semi-B-Fredholm operators) such that $\mbox{ind}(T)\leq 0$ (resp.
$\mbox{ind}(T)\geq 0).$ The spectrum associated to
$SBF_{+}^{-}(X)$ is called the {\it semi-essential approximate point
spectrum} and it is noted $ \sigma_{SBF_{+}^{-}}(T)=\{\lambda \in
\cit \,\, :\,\, T-\lambda \notin SBF_{+}^{-}(X) \}.$ While the
spectrum associated to $SBF_{-}^{+}(X)$ is noted $
\sigma_{SBF_{-}^{+}}(T)=\{\lambda \in \cit \,\, :\,\, T-\lambda
\notin SBF_{-}^{+}(X) \}.$

The  {\it ascent} $a(T)$ and the {\it descent}
$d(T)$ of $T$ are given by $a(T)=\inf\{n:N(T^{n})=N(T^{n+1})\},
d(T)=\inf\{n:R(T^n)=R(T^{n+1})\},$ with $\inf \emptyset = \infty.$ It is well-known
that if $a(T)$ and $d(T)$ are both finite then
they are equal, see \cite[Proposition 38.3]{Heuser}.

Recall that an operator $T$ is {\it Drazin invertible} if it
has a finite ascent and descent. It is well known that $T$ is
Drazin invertible if and only if  $T=R\oplus N$ where $R$ is invertible and $N$
is nilpotent (see \cite[Corollary 2.2]{Lay}). The Drazin spectrum is defined by
$\sigma_D(T)= \{ \lambda \in \mathbb{C}: T-\lambda I \,\,\, \text
{is not Drazin invertible} \}.$ From \cite[Lemma 4.1]{Be1} and \cite[Corollary 2.2]{Lay},
we have $$\sigma_{BW}(T)\subseteq\sigma_{D}(T).$$

Define the set $LD(X)$ as
follow : $$LD(X)=\{T\in {\mathcal B}(X) : a(T)<\infty \mbox{ and }
R(T^{a(T)+1}) \mbox{ is closed}\}.$$ From \cite{MM}, $LD(X)$ is a
regularity and it is the dual version of the regularity
$RD(X)=\{T\in {\mathcal B}(X) : d(T)<\infty \mbox{ and }
R(T^{d(T)}) \mbox{ is closed}\}.$ An operator $T\in {\mathcal B}(X)$ is said to be
{\it left} (resp. {\it right}) {\it
Drazin invertible} if $T\in LD(X)$ (resp. $T\in RD(X)$). The {\it left Drazin spectrum} $\sigma_{lD}(T)$ and the {\it right Drazin spectrum} $\sigma_{rD}(T)$ are defined by
$$\sigma_{lD}(T)=\{ \lambda\in \mathbb{C} \ :\ T-\lambda \notin
LD(X)\} \mbox{ and } \sigma_{rD}(T)=\{ \lambda\in \mathbb{C} \ :\
T-\lambda \notin RD(X)\}.$$ It is not difficult to see that $$\sigma_D(T)=\sigma_{lD}(T)\cup\sigma_{rD}(T).$$

\section{Preliminaries results}
An operator $T\in \mathcal {B}(X)$ has the single-valued
extension property at $\lambda_0\in \cit,$ SVEP at $\lambda_0,$ if
for every open disc $D_{\lambda_0}$ centered at $\lambda_0$ the
only analytic function $f\,\,:\,\, D_{\lambda_0}\longrightarrow X$
which satisfies $(T-\lambda)f(\lambda)=0$ for all $\lambda \in
D_{\lambda_0}$ is the function $f\equiv 0.$ Trivially, every
operator $T$ has SVEP at pint of the resolvent; also $T$ has SVEP
at $\lambda \in iso\,\sigma(T).$ We say that $T$ has SVEP if
it has SVEP at every $\lambda \in \cit$, \cite{Fi}. We denote by ${\mathcal S}(T)$ the
set of all $\lambda \in \cit$ such that $T$ does not have the
single-valued extension property at $\lambda.$ Note that (see \cite{Fi, LN})
\begin{equation}\label{surj}
{\mathcal S}(T)\subseteq\sigma_p(T) \mbox{ and } \sigma(T)={\mathcal S}(T)\cup\sigma_s(T).
\end{equation} In particular, if $T$ (resp. $T^*$) has the SVEP then $\sigma(T)=\sigma_s(T)$ (resp. $\sigma(T)=\sigma_a(T)$).

Since the ascent implies the SVEP (\cite{Laursen}) then we have
$${\mathcal S}(T)\subseteq\sigma_{lD}(T)\mbox{ and }{\mathcal S}(T^*)\subseteq\sigma_{rD}(T).$$

In the following theorem, we prove equality up to ${\mathcal S}(T)$
between the left Drazin spectrum and the left B-Fredholm spectrum.

\begin{thm}\label{t1} Let $T\in\mathcal{B}(X)$. Then
$$\sigma_{lD}(T)=\sigma_{lBF}(T)\cup {\mathcal S}(T).$$
\end{thm}
\begin{proof}
Let $\lambda\notin\sigma_{lD}(T)$, then
$R((T-\lambda)^{a(T-\lambda)+1})$ is closed. Hence
$R((T-\lambda)^{a(T-\lambda)})$ is closed by \cite[Lemma 12]{MM}.
Let $x\in N((T-\lambda)_{[a(T-\lambda)]})$ then $x\in
N(T-\lambda)\cap R((T-\lambda)^{a(T-\lambda)})$. Hence
$x=(T-\lambda)^{a(T-\lambda)}y$ for some $y\in X$. Hence
$0=(T-\lambda)x=(T-\lambda)^{a(T-\lambda)+1}y$. Thus $y\in
N((T-\lambda)^{a(T-\lambda)+1})=N((T-\lambda)^{a(T-\lambda)}).$
Therefore $x=0$ and hence $(T-\lambda)_{[a(T-\lambda)]}$ is
injective. On the other hand,
$R((T-\lambda)_{[a(T-\lambda)]})=R((T-\lambda)^{a(T-\lambda)+1})$
is closed. Thus $(T-\lambda)_{[a(T-\lambda)]}$ is upper
semi-Fredholm and hence $\lambda\notin \sigma_{lBF}(T).$ Since we
have ${\mathcal S}(T)\subseteq\sigma_{lD}(T)$ then
$$\sigma_{lBF}(T)\cup {\mathcal S}(T)\subseteq \sigma_{lD}(T).$$ Now let
$\lambda\notin[\sigma_{lBF}(T)\cup ({\mathcal S}(T)]$, then $T$ has the SVEP
at $\lambda$ and $T-\lambda$ is upper semi-B-Fredholm operator.
Hence it follows from \cite[Proposition 3.2]{Be2} that there exist $n$ such that
$R((T-\lambda)^n)$ is closed and $(T-\lambda)_{[n]}$ is semi regular.
Since $(T-\lambda)_{[n]}$ has also the SVEP at $0$, then from
\cite[Theorem 3.14]{Ai}, we conclude that $(T-\lambda)_{[n]}$ is injective with
closed range. Let $x\in N(T-\lambda)^{n+1}$, then
$(T-\lambda)(T-\lambda)^nx=0$. Hence $(T-\lambda)^nx\in
N(T-\lambda)\cap R(T-\lambda)^n=N((T-\lambda)_{[n]})=\{0\}.$ Thus
$x\in N(T-\lambda)^n$ , and hence
$N(T-\lambda)^{n}=N(T-\lambda)^{n+1}.$ So $T-\lambda$ is of finite
ascent and $a(T-\lambda)\leq n.$ We have
$R((T-\lambda))_{[n]}\subseteq N(T-\lambda)^{n+1}$ with
$a(T-\lambda)+1\leq n+1$. Therefore
$R(T-\lambda)^{a(T-\lambda)+1}$ is closed. Thus $T-\lambda$ is
left Drazin invertible. Hence
$\sigma_{lD}(T)\subseteq\sigma_{lBF}(T)\cup {\mathcal S}(T).$
\end{proof}

\begin{cor} If $T\in\mathcal{B}(X)$ has the SVEP then
$$\sigma_{lD}(T)=\sigma_{lBF}(T).$$
\end{cor}

By duality we get a similarly result for the right
Drazin spectrum.

\begin{thm}\label{t2} Let $T\in\mathcal{B}(X)$. Then
 $$\sigma_{rD}(T)=\sigma_{rBF}(T)\cup {\mathcal S}(T^*).$$
\end{thm}

\begin{proof}
Since $\sigma_{rBF}(T)=\sigma_{lBF}(T^*)$ and $\sigma_{rD}(T)=
\sigma_{lD}(T^*).$ We conclude by Theorem
 \ref{t1}.
\end{proof}

\begin{cor} If $T^*\in\mathcal{B}(X)$ has the SVEP then
$$\sigma_{rD}(T)=\sigma_{rBF}(T).$$
\end{cor}

From Theorem \ref{t1} and Theorem \ref{t2} we get the following
corollary.

\begin{cor}Let $T\in\mathcal{B}(X)$. Then
$$\sigma_{D}(T)=\sigma_{BF}(T)\cup [{\mathcal S}(T)\cup {\mathcal S}(T^*)].$$
\end{cor}

\begin{cor}\label{c20} If $T$ and $T^*$ have the SVEP then
$$\sigma_{D}(T)=\sigma_{BF}(T).$$
\end{cor}

In the following theorem, we prove equality up to ${\mathcal S}(T)$
between the left Drazin spectrum and the semi-essential
approximate point spectrum.

\begin{thm}\label{t3}Let $T\in\mathcal{B}(H)$ where $H$ is a Hilbert space. Then
 $$\sigma_{SBF_+^-}(T)\cup {\mathcal S}(T)=\sigma_{lD}(T).$$
\end{thm}

\begin{proof}From \cite[Lemma 2.12]{BK} we have
$\sigma_{SBF_+^-}(T)\subseteq\sigma_{lD}(T)$ and since ${\mathcal S}(T)
\subseteq \sigma_{lD}(T)$. Then
$$\sigma_{SBF_+^-}(T)\cup {\mathcal S}(T)\subseteq\sigma_{lD}(T).$$ Now let
$\lambda\in\sigma_{lD}(T)\setminus\sigma_{SBF_+^-}(T)$. We can
assume that $\lambda=0$, then $T$ is an upper semi-B-Fredholm,
$\mbox{ind}(T)\leq 0$ and $T$ is not left Drazin invertible. Hence
$a(T)=\infty$ or $R(T^{a(T)+1})$ is not closed. If $a(T)$ is
finite, then necessarily $R(T^{a(T)+1})$ is closed since $T$ is
upper semi-B-Fredholm. So $a(T)=\infty$. Also from \cite{BK}
$T=F\oplus N$ where $F$ is an upper semi-Fredholm operator and $N$
is a nilpotent operator. Since $N$ is nilpotent and $a(T)=\infty$,
then $a(F)=\infty$. Now it follows from \cite{Ai} that $F$ does not
have the SVEP at $0$, and then nor $T$.
\end{proof}

\begin{cor} If $T\in\mathcal{B}(H)$ has the SVEP then
 $$\sigma_{SBF_+^-}(T)=\sigma_{lD}(T).$$
\end{cor}

By duality we have the following result :

\begin{thm}\label{t4} Let $T\in\mathcal{B}(H)$.
Then $$\sigma_{SBF_-^+}(T)\cup {\mathcal S}(T^*)=\sigma_{rD}(T).$$
\end{thm}
\begin{proof}
Since $\sigma_{SBF_+^-}(T^*)=\sigma_{SBF_-^+}(T)$ and
$\sigma_{lD}(T^*) =\sigma_{rD}(T).$ We conclude by Theorem
\ref{t3}.\end{proof}

\begin{cor} If $T^*\in\mathcal{B}(H)$ has the SVEP then
$$\sigma_{SBF_-^+}(T)=\sigma_{rD}(T).$$
\end{cor}

From Theorem \ref{t3} and Theorem \ref{t4},  it is not hard to see that
$$\sigma_{D}(T)=\sigma_{BW}(T)\cup [{\mathcal S}(T)\cup {\mathcal S}(T^*)].$$
In fact, more can be said in the Banach setting
\begin{thm} \label{c-2}Let $T\in\mathcal{B}(X)$ then
$$\sigma_{D}(T)=\sigma_{BW}(T)\cup [{\mathcal S}(T)\cap {\mathcal S}(T^*)].$$
\end{thm}
\begin{proof}Since $\sigma_{BW}(T)\cup ({\mathcal S}(T)\cap {\mathcal S}(T^*))\subseteq\sigma_D(T)$ always holds,
then let $\lambda\notin\sigma_{BW}(T)\cup ({\mathcal S}(T)\cap {\mathcal S}(T^*)).$ Without loss of generality we assume that $\lambda=0$.
Then $T$ is a B-Fredholm operator of index zero.

{\it Case.1} If $0\notin {\mathcal S}(T)$ : Since $T$ is a
B-Fredholm operator of index zero, then it follows from
\cite[Lemma 4.1]{Be1} that there exists a Fredholm operator $F$ of
index zero and a nilpotent operator $N$ such that $T=F\oplus N$.
If $0\notin\sigma(F)$, then $F$ is invertible and hence $T$ is
Drazin invertible. Now assume that $0\in\sigma(F)$. Since $T$ has
the SVEP at $0$, then $F$ has also the SVEP at $0$. Hence it
follows from \cite[Theorem 3.16]{Ai} that $a(F)$ is finite. $F$ is
a Fredholm operator of index zero, then it follows from
\cite[Theorem 3.4]{Ai} that $a(F)$ is also finite. Then
$a(F)=d(F)<\infty$ which implies that $0$ is a pole of $F$ and
hence an isolated point of $\sigma(F)$. $N$ is nilpotent, then $0$
is isolated point of $\sigma(T)$. From \cite[Theorem 4.2]{Be1} we
get $0\notin\sigma_D(T)$. {\it Case.2} If $0\notin {\mathcal
S}(T^*)$, then proof goes similarly.
\end{proof}

\begin{cor}\label{dbwsvep}{\rm \cite{BCD}} If $T$ or $T^*$ has the SVEP then $$\sigma_{D}(T)=\sigma_{BW}(T).$$
\end{cor}

Recall that $T$ is a {\it Browder} operator if $T$ is a Fredholm
of finite ascent and descent. Let $\sigma_B(T)$ be the {\it
Browder spectrum} defined as the set of all $\lambda\in\mathbb{C}$
such that $T-\lambda$ is not Browder. Analogously, $T$ is {\it
B-Browder} operator if for some integer $n$, $R(T-\lambda)^n$ is
closed and $(T-\lambda)_{[n]}$ is Browder and let $\sigma_{BB}(T)$
be the {\it B-Browder spectrum}. In \cite[Corollary 3.53]{Ai} it
is proved that if $T$ or $T^*$ has the SVEP, then
$$\sigma_W(T)=\sigma_B(T).$$ From \cite[Theorem 3.6]{Be2}
$\sigma_D(T)=\sigma_{BB}(T)$, then by Corollary \ref{dbwsvep}, if
$T$ or $T^*$ has the SVEP then $$\sigma_{BW}(T)=\sigma_{BB}(T).$$
\begin{thm}\label{smt}Let $T\in\mathcal{B}(X)$ and $f$ be an analytic
 function on some open neighborhood of
$\sigma(T)$ which is nonconstant on any connected component of $\sigma(T)$ then
$$\sigma_{BW}(f(T))\cup [{\mathcal S}(f(T))\cap {\mathcal S}(f(T^*))]=f(\sigma_{BW}(T)\cup [{\mathcal S}(T)\cap {\mathcal S}(T^*)])$$
\end{thm}
\begin{proof} Since the Drazin spectrum satisfies the spectral mapping theorem for every analytic
 function $f$ on some open neighborhood of
$\sigma(T)$ which is nonconstant on any connected component of $\sigma(T)$,
then the result follows at once from Theorem \ref{c-2}.
\end{proof}

It is well known that if $T$ has the SVEP then $f(T)$ has also the
SVEP \cite{LN}. Now we retrieve the result proved in \cite{AM1,
HZ} : $f(\sigma_{BW}(T))=\sigma_{BW}(f(T))$ whenever $T$ or $T^*$
has the SVEP. Note that in \cite{AM1, HZ} the condition " $f$ is
nonconstant on any connected component of $\sigma(T)$" is dropped.

\section{Applications }

\subsection{Perturbations}

\begin{lem} \label{l1}Let $T\in\mathcal{B}(X)$. Let $N\in\mathcal{B}(X)$ be a nilpotent operator such that $TN=NT$.
Then $$\mathcal{S}(T+N)=\mathcal{S}(T).$$
\end{lem}
\begin{proof}See for instance \cite[Lemma 2.1]{BZZ2}
\end{proof}

\begin{lem}\label{l2} Let $T\in\mathcal{B}(X)$.
If $N\in\mathcal{B}(X)$ is a nilpotent operator which commutes
with $T$.
 Then $$\sigma_{lD}(T+N)=\sigma_{lD}(T).$$
\end{lem}

\begin{proof} Assume that $\lambda=0\notin \sigma_{LD}(T).$
Then $a(T)$ is finite and $R(T^{a(T)+1})$ is closed. Let $m$ be
the nonnegative integer such that $N^m=0\neq N^{m-1}.$ Let
$s=\mbox{max}(a(T),m)$. Then
\[
\begin{array}{ccl}
(T+N)^{2s}&=&{\displaystyle\Sigma_{k=0}^{2s}}(_k^{2s})T^kN^{2s-k}\\
&=&(_0^{2s})N^{2s}+\cdots+(_s^{2s})T^sN^s+(_{s+1}^{2s})T^{s+1}N^{s-1}+\cdots+(_{2s}^{2s})T^{2s}\\
&=&(_{s+1}^{2s})T^{s+1}N^{s-1}+\cdots+(_{2s}^{2s})T^{2s}\\
&=&T^s²[(_{s+1}^{2s})T^{1}N^{s-1}+\cdots+(_{2s}^{2s})T^{s}].\\
\end{array}\]
Now Let $x\in N(T)^{2s}=N(T)^{s}$ that is $(T)^{2s}x=0$. Then it
follows form the above equality that $(T+N)^{2s}x=0$. Hence
$N(T)^{2s}\subseteq N(T+N)^{2s}.$ With the same argument for $T+N$
and $-N$ we have $N(T+N)^{2s}\subseteq N(T)^{2s}.$ Thus
$N(T)^{2s}=N(T+N)^{2s}.$ Since $N(T^s)=N(T^{2s})=N(T^{2s+1})$,we
get $N(T+N)^{2s}=N(T+N)^{2s+1}.$ Therefore $T+N$ is of finite
ascent. In the other hand $R(T+N)^{2s}\subseteq R(T^s)$ is closed.
Hence by \cite[Lemma 12]{MM} $R(T+N)^{2s+1}$ is closed. Thus
$0\notin \sigma_{lD}(T+N).$
\end{proof}
The following result follows from Theorem \ref{t3}, Lemma \ref{l1} and
Lemma \ref{l2}
\begin{thm}Let $T\in\mathcal{B}(X)$. Let $N\in\mathcal{B}(X)$
be a nilpotent operator which commutes with $T$. Then
$$\sigma_{SBF_+^-}(T+N)\cup\mathcal{S}(T)=\sigma_{SBF_+^-}(T)\cup\mathcal{S}(T).$$
\end{thm}

The following corollary which is proved in \cite{AM4} gives an
affirmative answer for the question posed by Berkani-Amouch
\cite{BA} in the case where $T$ has the SVEP.

\begin{cor}Let $T\in\mathcal{B}(X)$ have the SVEP. Let $N\in\mathcal{B}(X)$
be a nilpotent operator which commutes with $T$. Then $$\sigma_{SBF_+^-}(T+N)=\sigma_{SBF_+^-}(T).$$
\end{cor}

\subsection{Generalized Weyl's theorem for operator matrices}

Berkani \cite[Theorem 4.5]{Be1} has shown that every
normal operator $T$ acting on Hilbert space $H$ satisfies
\begin{equation}\label{eqgw}\sigma(T)\setminus
E(T)=\sigma_{BW}(T),\end{equation}where $E(T)$ is the set of all
isolated eigenvalues of $T$. We say that {\it generalized Weyl's
theorem} holds for $T$ if Equation (\ref{eqgw}) holds. This gives
a generalization of the classical Weyl's theorem. Recall that $T
\in {\mathcal B}(X)$ obeys {\it Weyl's theorem} if
\begin{equation}\label{weyy}
\sigma(T)\setminus E_0(T)=\sigma_{W}(T),
\end{equation} where
$E_0(T)$ denotes the set of the isolated points of $\sigma(T)$
which are eigenvalues of finite multiplicity. Form \cite[Theorem
3.9]{BK} generalized Weyl's theorem implies Weyl's theorem and
generally the reverse is not true.

For $A\in{\mathcal B}(X),\,B\in{\mathcal B}(Y)$ and $C\in{\mathcal B}(Y,X)$ we
 denote by $M_C$ the operator defined
 on $X\oplus Y$ by
 $$M_C=\begin{bmatrix}
   A & C \\
 0 & B
 \end{bmatrix}.$$
In general the fact that
generalized Weyl's theorem holds for $A$ and $B$ does not imply
that generalized Weyl's theorem holds for
$M_0=[\begin{smallmatrix}
   A & 0 \\
 0 & B
 \end{smallmatrix}]$. Indeed, Let $I_1$ and $I_2$ be the identities on $\cit$ and
$l_2$, respectively. Let $S_1$ and $S_2$ defined on $l_2$ by
$$S_1(x_1,x_2,\ldots)=(0,{1\over 3}x_1,{1\over 3}x_2,\ldots),\quad
S_2(x_1,x_2,\ldots)=(0,{1\over 2}x_1,{1\over 3}x_2,\ldots).$$Let
$T_1=I_1\oplus S_1,$ $T_2=S_2-I_2,$ $A=T_1^2$ and $B=T_2^2$, then
from \cite[Example 1]{HZ} we have $A$ and $B$ obey generalized
Weyl's theorem but $M_0$ does not obey it. It also may happen that
$M_C$ obeys generalized Weyl's theorem while $M_0$ does not obey
it. Let $A$ be the unilateral unweighed shift operator. For
$B=A^*$ and $ C = I- AA^*$, we have that $M_C$ is unitary without
eigenvalues. Hence $M_C$ satisfies generalized Weyl's theorem (see
\cite[Remark 3.5]{Berkani2004}). But
$\sigma_{w}(M_0)=\{\lambda\,:\,|\lambda|=1\}$ and
$\sigma(M_0)\setminus E_0(M_0)=\{\lambda\,:\,|\lambda|\leq 1\}.$
Then $M_0$ does not satisfy Weyl's theorem and so from
\cite[Theorem 3.9]{BK} it does not satisfy generalized Weyl's
theorem either.

A bounded linear operator $T$ is said to be {\it
isoloid} if every isolated point of $\sigma(T)$ is an eigenvalue
of $T$.

\begin{prop}\label{prop3.1} Let $A$ and $B$ be isoloid.
Assume that $\sigma_{BW}(A\oplus B)=\sigma_{BW}(A)\cup\sigma_{BW}(B)$.
If $A$ and $B$ obeys generalized Weyl's theorem, then $A\oplus B$ obeys generalized Weyl's theorem.
\end{prop}
\begin{proof}Since $A$ and $B$ are isoloid, then
$$E(A\oplus B)=[E(A)\cap \rho(B)]\cup[\rho(A)\cap E(B)]\cup[E(A)\cap E(B)].$$
Now if $A$ and $B$ obeys generalized Weyl's theorem, then
\[
\begin{array}{ccl}
E(A\oplus B)&=&[\sigma(A)\cup\sigma(B)]\setminus[\sigma_{BW}(A)\cup\sigma_{BW}(B)]\\
&=&\sigma(A\oplus B)\setminus\sigma_{BW}(A\oplus B).\\
\end{array}\]
Then $A\oplus B$ obeys generalized Weyl's theorem.
\end{proof}

\begin{lem}\label{BWMatrice}Let $A\in{\mathcal B}(X)$ and $B\in{\mathcal B}(Y)$ have the SVEP, then $$\sigma_{BW}(M_C)=\sigma_{BW}(A)\cup\sigma_{BW}(B),$$ for all $C\in{\mathcal B}(Y,X).$
\end{lem}
\begin{proof}Since $A$ and $B$ have the SVEP, then it follows from \cite[Proposition 3.1]{Zg2} that $M_C$ has also the SVEP. Hence $\sigma_{BW}(M_C)=\sigma_{D}(M_C)$ by Corollary \ref{dbwsvep}. Also since $A$ and $B$ have the SVEP, it follows from \cite[Corollary 2.1]{HZ2} that $\sigma_{D}(M_C)=\sigma_{D}(A)\cup\sigma_{D}(B)$. Therefore $\sigma_{BW}(M_C)=\sigma_{BW}(A)\cup\sigma_{BW}(B)$ by Corollary \ref{dbwsvep}.
\end{proof}

\begin{thm}\label{t-1} Let $A$ and $B$ be isoloid with the SVEP.
If $A$ and $B$ obeys generalized Weyl's theorem,
then $M_C$ obeys generalized Weyl's theorem for every $C\in{\mathcal B}(Y,X).$
\end{thm}
\begin{proof}It follows from Proposition \ref{prop3.1} and Lemma \ref{BWMatrice} that
$$E(A\oplus B)=\sigma(A\oplus B)\setminus\sigma_{BW}(A\oplus B)=\sigma(M_C)\setminus\sigma_{BW}(M_C).$$
So it is enough to show that $E(A\oplus B)=E(M_C).$ Let $\lambda\in E(M_C)$.
Then $\lambda\in\sigma_p(M_C)\subseteq\sigma_p(A)\cup\sigma_p(B).$
Hence $\lambda\in\sigma_p(A\oplus B).$ Since $\lambda\in iso\,\sigma(M_C)=iso\,\sigma(A\oplus B)$
then $\lambda\in E(A\oplus B).$ Now Let $\lambda\in E(A\oplus B).$ If $\lambda\in\sigma(A)$
then $\lambda\in iso\,\sigma(A).$ Since $A$ is isoloid, then $\lambda\in\sigma_p(A)\subseteq\sigma_p(M_C).$
Hence $\lambda\in E(M_C).$ If $\lambda\in\sigma(B)\setminus\sigma(A)$, then $\lambda\in\sigma_p(B).$
Since $A$ is invertible, then $\lambda\in\sigma_p(M_C).$ Thus $\lambda\in E(M_C).$
 Therefore $E(A\oplus B)=E(M_C).$
\end{proof}

Let $\pi(T)$ be the set of all poles of the resolvent of $T$. Recall from \cite{DHJ} that $T$ is {\it polaroid} if  $iso\,\sigma(T)\subseteq \pi(T)$. Since $\pi(T)\subseteq E(T)$ holds without restriction on $T$, then if $T$ is polaroid then $E(T)=\pi(T)$.

\begin{cor}Let $A$ and $B$ be polaroid with the SVEP.
Then $M_C$ obeys generalized Weyl's theorem for every
$C\in{\mathcal B}(Y,X).$
\end{cor}
\begin{proof}
$A$ and $B$ are polaroid then $E(A)=\pi(A)$ and $E(B)=\pi(B)$. Since $A$ and $B$ has the SVEP, then by
\cite{AZ}, $A$ and $B$ satisfies generalized Weyl's
theorem. Hence we conclude by Theorem \ref{t-1}.
\end{proof}

\subsection{Multipliers on a commutative Banach algebra}
Let $\mathcal{A}$ be a semi-simple commutative Banach algebra. A mapping $T\,:\,\mathcal{A}\longrightarrow \mathcal{A}$ is called a {\it multiplier} if $$T(x)y=xT(y)\mbox{ for all }x,\,y\in \mathcal{A}.$$ By semi-simplicity of $\mathcal{A}$, every multiplier is a bounded linear operators on $\mathcal{A}$. Also the semi-simplicity of $\mathcal{A}$ implies that every multiplier has the SVEP (see \cite{Ai, LN}).

From Corollary \ref{dbwsvep} we have
\begin{prop}{\rm\cite{BA2}} Let $T$ be a multiplier on semi-simple commutative Banach algebra $\mathcal{A}$, then the following assertions are equivalent
\begin{itemize}
\item[i)] $T$ is B-Fredholm of index zero.
\item[ii)] $T$ is Drazin invertible.
\end{itemize}
\end{prop}
From \cite[Theorem 4.36]{Ai}, for every multiplier $T$ on
semi-simple commutative Banach algebra $\mathcal{A}$,
$E(T)=\pi(T)$ and since $T$ has the SVEP we get from \cite{AZ}
\begin{prop}Every multiplier on semi-simple commutative Banach algebra $\mathcal{A}$ obeys generalized Weyl's thoerem.
\end{prop}
Also from Theorem \ref{smt}
\begin{cor}{\rm\cite{BA2}} Let $T$ be a multiplier on semi-simple commutative Banach algebra $\mathcal{A}$. Let $f$ be an analytic
 function on some open neighborhood of
$\sigma(T)$ which is nonconstant on any connected component of $\sigma(T)$ then
$$\sigma_{BW}(f(T))=f(\sigma_{BW}(T)).$$
\end{cor}
Now if assume in additional that $\mathcal{A}$ is regular and Tauberian (see \cite{LN} for definition) then every multiplier $T$ has the weak decomposition property ($\delta_w$) and then $T^*$ has also the SVEP (see \cite{ZZ2005} for definition and details). Hence we get from Corollary \ref{c20}
\begin{prop}Let $T$ be a multiplier on semi-simple regular Tauberian commutative Banach algebra $\mathcal{A}$, then the following assertions are equivalent
\begin{itemize}
\item[i)] $T$ is B-Fredholm operator.
\item[ii)] $T$ is Drazin invertible.
\end{itemize}
\end{prop}
For $G$ a locally compact abelian group, let $L^1(G)$ be the space
of $\cit$-valued functions on $G$ integrable with respect to Haar
measure and $M(G)$ the Banach algebra of regular complex Borel
measures on $G$. We recall that $L^1(G)$ is a regular semi-simple
 Tauberian commutative Banach algebra. Then we have the following
\begin{cor}Let $G$ be a locally compact abelian group, $\mu\in M(G)$ and
$X=L^1(G)$. Then every  convolution operator
$T_{\mu}\,:\,X\longrightarrow X$, $T_{\mu}(k)=\mu\star k$ is B-Fredholm if and only if is Drazin invertible.
\end{cor}

\end{document}